\documentclass[11pt]{amsart}
\usepackage{amsmath}
\usepackage{amssymb}
\usepackage{enumerate}

\newtheorem{thm}{Theorem}
\newtheorem{lem}[thm]{Lemma}
\newtheorem{prop}[thm]{Proposition}

\begin{document}

\title{Birkhoff spectra for one-dimensional maps 
with some hyperbolicity}

\author{Yong Moo Chung }

\address{ Department  of  Applied  Mathematics, 
Hiroshima  University,
Higashi-Hiroshima 739-8527, Japan}

\curraddr{Department of Mathematics,
The Pennsylvania State University \\
University Park, PA 16802, USA}

\email{chung@amath.hiroshima-u.ac.jp}

\date{February 5, 2008. 
{\it Last modified:} February 25, 2008.
\\
\hspace{0.4cm}{\it 2000 Mathematics Subject Classification.} 
Primary 37C45; Secondary 37D25, 37E05.}



\begin{abstract}
We study the multifractal analysis for
smooth dynamical systems in dimension one.
It is characterized 
the Hausdorff dimension of the level set obtained from
the Birkhoff averages of a continuous function 
by the local dimensions of hyperbolic measures
for a topologically mixing $C^2$ map modelled by 
an abstract dynamical system.
A characterization which corresponds to above is also given 
for the ergodic basins of invariant probability measures.
And it is shown that 
the complement of the set 
of quasi-regular points
carries full Hausdorff dimension. 
\end{abstract}

\pagestyle{plain} 

\maketitle

\section{Introduction}

Multifractal analysis is a theory to understand  
phenomena associated with more than one scaling symmetry,
and a practical tool
for the numerical study of dynamical systems.
The purpose of mathematical study of multifractal analysis
is to investigate the multifractal spectra,
that is the quantities concerning with the level sets
obtained from invariant local quantities. 
We refer to \cite{P2} for a general account of the subject.
The multifractal spectra is well understood 
for conformal (or low dimensional) hyperbolic dynamical systems
\cite{Ol, PW97, PW}.  
Also, some of solid mathematical results  
are known for some class of smooth one-dimensinal maps 
with indifferent periodic points
such as the Manneville-Pomeau maps \cite{I, Na, PoW}.
Takens and Verbitskiy gave a
variational principle concerning with
the entropy spectra for Birkhoff averages
by assuming the specification property for dynamical systems \cite{TV2}.   
However, it seems that there are no rigorous results
on the dimension spectra
for smooth maps with critical points
such as logistic maps.
This paper is devoted to study the dimension spectra
for the Birkhoff averages of 
smooth interval maps without assuming uniform hyperbolicity.

Let $I$ be a compact interval of the real line $\mbox{\boldmath $R$}$
and $f:I\to I$ a $C^2$ map.
For a repeller $\varLambda$ for $f$,
it is known that 
$\dim_H (\varLambda )$, the Hausdorff dimension of $\varLambda$,
is given by the unique solution $\delta$ of the Bowen equation
$P_{\varLambda} (\delta\cdot \log |f'| )=0,$
where $P_{\Lambda}$ denotes the topological pressure
for the restriction of $f$ to $\varLambda$ \cite{B1},
and then it is obtained a dimension formula
\begin{equation}
\label{bowen:f} \dim_H (\varLambda ) 
=\sup \left\{
\frac{h_{\mu} (f)}{\int \log |f'| d\mu} :
\mu\in\mathcal{M}_f , \mu (\varLambda)=1 \right\}
\end{equation}
by the variational principle \cite{W},
where $\mathcal{M}_f$ denotes 
the set of the $f$-invariant Borel probability measures 
on $I$
and $h_{\mu} (f)$ the metric entropy of $\mu\in\mathcal{M}_f .$
Notice that the supremum is attained by an ergodic measure.
In this paper we give a formula which corresponds to 
(\ref{bowen:f}) for the level sets obtained by the Birkhoff averages
of continuous functions for 
a class of smooth interval maps modelled by abstract dynamical systems.

We need some notions and facts to state our results precisely.  
For a continuous function $\varphi : I\to  \mbox{\boldmath $R$}$
we denote by $B_{\varphi} (\alpha)$ the set of points $x\in I$
such that $S_n\varphi (x) / n$ converges to $\alpha\in\mbox{\boldmath $R$}$
as $n\to\infty$,
where $$S_n\varphi (x) :=\varphi (x) + \varphi (f(x)) +
\cdots + \varphi (f^{n-1}(x)) .$$
Then the phase space $I$ is decomposed into the level sets 
$B_{\varphi} (\alpha)$
as follows:
$$I= 
\left( \bigsqcup_{\alpha\in\mbox{\boldmath $R$}} B_{\varphi} (\alpha) \right)
\sqcup I_{\varphi}$$
where $I_{\varphi}$ denotes the irregular set for $\varphi$, 
that is the set of points $x\in I$ such that
$S_n\varphi (x) / n$
does not converge as $n\to \infty$.
We say that $\mu\in \mathcal{M}_f$ is {\it hyperbolic} if
the Lyapunov exponent
$\displaystyle \limsup_{n\to\infty} \frac{1}{n}\log |(f^n)'(x) | $
is positive for $\mu$-almost every $x\in I .$
If $\mu\in \mathcal{M}_f$ is ergodic 
then the Lyapunov exponent is constant $\mu$-almost everywhere,
and the constant is given by
$$\displaystyle \lambda_{\mu} (f):= \int \log |f'| d\mu .$$ 
We call $\lambda_{\mu} (f)$ {\it the Lyapunov exponent} of
$\mu$,
and denote by $\mathcal{H}_f$ the set of 
$f$-invariant ergodic probability measures $\mu$ which is hyperbolic,
i.e. $\lambda_{\mu} (f)>0 .$
The Ruelle inequality \cite{Ru} asserts that
$$\lambda_{\mu} (f)\ge h_{\mu} (f)$$
holds whenever $h_{\mu} (f)>0 ,$
and hence any ergodic measure with positive metric entropy is hyperbolic.
Moreover, a theory for hyperbolic measures \cite{C0,  K, KM} 
asserts the following:
\begin{prop}
Let $\mu\in\mathcal{H}_f .$ 
Then for any
$\varphi_1 ,\ldots , \varphi_p \in C(I)$
and $\varepsilon >0$ 
there are an integer $k\ge 1$,
a compact interval $L\subset I$ and
a family $\mathcal{K}$ of pairwise disjoint compact intervals
with $K\subset L =f^k(K)$ for each $K\in\mathcal{K}$ 
such that
$$(\log \sharp\mathcal{K} )/k\ge h_{\mu} (f) -\varepsilon ,
\quad
| \frac{1}{k} \log\vert (f^k)' (x) \vert 
-\lambda_{\mu} (f) | \le\varepsilon $$
and  
$$|\frac{1}{k}S_k \varphi_j (x) - \int \varphi d\mu |\le
\varepsilon \enspace (j=1,\ldots , p)$$
for all $x\in\sqcup_{K\in\mathcal{K}} K,$
where 
$C(I)$ denotes the space of the continuous functions on $I$,
and $\sharp A$ the cardinality of a set $A$.  
\end{prop}
For $\mu\in\mathcal{H}_f$
we call the quantity $D(\mu)$ given by
$$D(\mu ):=\frac{h_{\mu} (f)}{\lambda_{\mu} (f)}$$ 
{\it the local dimension of $\mu$},
because 
$$ D(\mu ) =\lim_{r\to 0}\frac{ \log \mu ([x-r , x+r])}{\log r}$$
holds for $\mu$-almost every $x\in I$
\cite{H, L2}.
Then it follows that
$$\dim_H (\mu ) = D(\mu )$$ 
where $\dim_H (\mu)$ denotes {\it the Hausdorff dimension of $\mu$}
defined by 
 $$\dim_H (\mu ):=\inf\{ \dim_H (Y) : Y\subset I 
 \text{ is a Borel set with }  \mu (Y)=1 \} $$ 
\cite{Y0}.
Notice that $D(\mu )=1$ holds if $\mu\in\mathcal{H}_f$
is absolutely continuous with respect to Lebesgue measure.
We say that a map $f:I\to I$ is {\it topologically mixing} if
for any nontrivial interval $L\subset I,$ 
$f^k(L)=I$ holds
whenever $k\ge 1$ is large.   
Throughout this paper we assume that $f:I\to I$ is topologically mixing.
Then the following lemma concerning with the structure of the set
of hyperbolic measures is established:
\begin{lem}
Let $f:I\to I$ be a topologically mixing $C^2$ map. 
Then ergodic measures are dense in the set of hyperbolic measures
with positive metric entropies,
i.e.,
for any hyperbolic $\mu\in\mathcal{M}_f$ with $h_{\mu } (f)>0,$  
$\varphi_1 ,\ldots , \varphi_p \in C(I)$
and small $\varepsilon >0$
there exists $\nu\in\mathcal{H}_f$ such that
$$|\int\varphi_j d\nu  - \int\varphi_j d\mu |\le \varepsilon$$
for each $j=1,2,\ldots , p.$
Moreover,  the measure $\nu$ can be taken to satisfy
$$h_{\nu} (f)\ge  h_{\mu} (f) - \varepsilon  \quad
\text{ and }\quad
\lambda_{\nu} (f) \le \int\log |f' |d\mu +\varepsilon .$$
\end{lem}
We need more assumptions on the map $f$ 
concerning with both of recurrence times and nonuniform hyperbolicity
to establish a dimension formula 
for the level sets obtained by Birkhoff averages.
We say that a function 
$R : J\to\mbox{\boldmath $N$}\cup \{ \infty \}$ is 
{\it a return time function} with the {base set} $J\subset I$
if $f^{R(x)}(x)\in J$ holds whenever $x\in J$ and $R(x)<\infty.$ 
Here, the function $R$ is not necessary 
to be the first return time to $J$.
We assume that the map $f: I\to I$ has a return time function
$R : J\to\mbox{\boldmath $N$}\cup \{ \infty \}$
with a compact interval $J\subset I$
as the base set satisfying the following
properties:
\begin{enumerate}[(H1)]
\item
there is a constant $\lambda >1$ such that for any $n\ge 1$
if $V$ is a connected component of $\{ x\in J: R(x) =n \}$, then 
 $$f^nV=J \quad \text{and} \quad 
|(f^n)'(x)|\ge \lambda \enspace (x\in V);$$ 
\item 
there is a sequence  
$\{\varepsilon_k\}_{k=0}^{\infty}$ of positive numbers which converges to
 zero as $k\to\infty$ such that for any $n\ge 1$
if 
$A$ is a connected component of
\begin{align*}
\{ x \in J : &R(x)= k_1 ,\enspace R(f^{k_1}(x))=k_2, \enspace 
\ldots \enspace ,\\
&R(f^{k_1+\cdots +k_{l-2}}(x))= k_{l-1},\enspace   
R(f^{k_1+\cdots +k_{l-1}}(x))\ge k_l  
\}
\end{align*} 
for some integers $k_1 , \ldots ,  k_l\ge 1$  
with $k_1+\cdots + k_l=n$,
then 
$$| f^i(A)| \le \varepsilon_{n-i}$$
holds for all $0\le i\le n-1$,
where $|A|$ denotes the length of an interval $|A|$;
\item 
there are a constant $C\ge 1$ and
a sequence  
$\{\eta_k\}_{k=0}^{\infty}$ of positive numbers which converges to
 zero as $k\to\infty$
such that for any $n\ge 1$ if 
$B$ is a connected component of
$$\{ x \in J :R(x)= k_1 ,\enspace R(f^{k_1}(x))=k_2, \enspace \ldots \enspace ,
\enspace R(f^{k_1+\cdots +k_{l-1}}(x))= k_l \} $$
for some integers $k_1 , \ldots ,  k_l\ge 1$  
with $k_1+\cdots + k_l=n ,$
then 
\begin{align*}
\frac{|(f^{n})'(y)|}{|(f^{n})'(z)|} \le C
\quad{ and }\quad
\frac{|(f^{k_1})'(y)|}{|(f^{k_1})'(z)|} \le e^{\eta_l}
\end{align*}
hold for all $y, z\in B$;  
\item 
there are an integer $l_0\ge 1$ and a constant $\gamma_0 >0$
such that for any $n\ge 1$
if $U$ is a connected component of $\{ x\in J: R(x) >n \}$, 
then
$$
 m \left( \left\{ x\in U : R(x)\le n+l_0 \right\} \right) \ge
\gamma_0 m(U) .$$
\end{enumerate}
A Collet-Eckmann unimodal map is an example 
which has a return time function 
satisfying the assumptions mentioned above 
\cite{Y3}.
A smooth interval map satisfying these assumptions 
is modelled by an abstract dynamical system on a tower 
having bounded slope \cite{C4}.
Then, the tail of the tower decays exponentially fast,
and then, an absolutely continuous invariant probability
measure exists and it 
satisfies the exponential decay of correlations \cite{Y3}. 
Moreover, the abstract dynamical system 
obtained from these assumptions has 
a kind of specification property, and from which
it is shown the large deviation principle \cite{C4}. 
The main result of this paper is the following:
\begin{thm}[the dimension spectra for Birkhoff averages]
Let $f: I\to I$ be a topologically mixing $C^2$ map
and assume that there is a return time function 
$R : J\to\mbox{\boldmath $N$}\cup \{ \infty \}$
satisfying (H1)-(H4).
Then for any 
$\varphi\in C(I)$
$$\dim_H (G_{\varphi} (\alpha)) 
= \lim_{\varepsilon\to 0+}
\sup{}^+ \Big\{ D(\nu) : \nu\in \mathcal{H}_f  , 
 |\int\varphi d\nu -\alpha | <\varepsilon \Big\} $$
holds for all $\alpha \in \mbox{\boldmath $R$} ,$ 
where $\sup^+ A := \max\{ \sup A , 0 \}$ for a set 
$A\subset\mbox{\boldmath $R$}$.
\end{thm}

\noindent {\bf Remark.}
The upper regularization in the right-hand side
of the formula above is not removable 
without assuming uniform hyperbolicity,
because the local dimension $D(\nu )$ may not be upper semicontinuous 
as a function of $\nu\in\mathcal{H}_f$.
In fact,
the Lyapunov exponent $\lambda_{\nu} (f)$ 
is not lower semicontinuous for $\nu$
if $f:I\to I$ is a Collet-Eckmann unimodal map
\cite{BK}.

\vspace{0.2cm}

It is also considerable the dimension spectra for 
the distributions along the orbits of dynamical systems. 
A point $x\in I$ is called \textit{quasi-regular for $f$} 
if the sequence of 
the empirical distributions 
$$\delta_{x}^n :=
\frac{1}{n} (\delta_x +\delta_{f(x)}+\cdots +\delta_{f^{n-1}(x)})
\in\mathcal{M}$$
converges to an invariant probability measure in the weak* toplogy
as $n\to\infty$,
where $\mathcal{M}$ denotes the space of the Borel probability measures
on $I$, and
$\delta_y\in\mathcal{M}$ the Dirac measure supported on $y\in I .$  
We denote by $QR(f)$ 
the set of the quasi-regular points for $f$.
For $\mu\in\mathcal{M}_f$ we call
the set 
$$B(\mu ) :=\{ x\in I : 
\lim_{n\to\infty} \delta_{x}^n =\mu \}$$
\textit{the ergodic basin of $\mu$}.
The set of the quasi-regular points is decomposed into 
the ergodic basins of invariant probability measures as follows:
$$QR(f) =\bigsqcup_{\mu\in\mathcal{M}_f}  B(\mu)  .$$
The next result of this paper asserts that 
the Hausdorff dimension of the ergodic basin 
for $\mu\in\mathcal{M}_f$ is represented 
by the local dimension of the hyperbolic ergodic measures
close to $\mu$.
It is the following:
\begin{thm}[the Hausdorff dimension of the ergodic basin]
Let $f:I\to I$ be as in Theorem 3. Then
for any $\mu\in\mathcal{M}_f$
$$\dim_H (B ( \mu ))  
=\inf_{\mathcal{U}}\sup\{
D(\nu ) : \nu\in\mathcal{H}_f\cap \mathcal{U} \}  $$
holds, where the infimum is taken over
all of the neighborhoods $\mathcal{U}$ of $\mu$ in $\mathcal{M} .$ 
\end{thm}

\noindent {\bf Remark.}
Neither the ergodicity nor the hyperbolicity is assumed
for $\mu\in\mathcal{M}_f$ in the theorem above.
Thus the local dimension  may not be defined 
for $\mu$ itself.

\vspace{0.2cm}

It is well-known that if $f:I\to I$ is topologically mixing, 
then the complement of $QR(f)$ belongs to the first category and
it is a null set for any invariant probability measure.
However, it is not a small set from the view point 
of dimension theory.
Another result of this paper corresponds to that 
obtained for uniformly expanding conformal maps 
by Barreira and Schmeling \cite{BS}:

\begin{thm}[the Hausdorff dimension of the irregular set]
Let $f: I\to I$ be as in Theorem 3. 
Then $$\dim_H ( I \setminus QR(f) ) =1 $$ holds.
\end{thm}


\noindent {\bf Acknowledgement.}
This paper was written while the author was 
visiting to the Pennsylvania State University
by a support of the Ministry of Education, 
Culture, Sports, Science and Technology of Japan.
He is grateful to Professors Ya. Pesin and O. Sarig for their
hospitality and for helpful conversertions. 
He also thanks to Professor N. Sumi for many fruitful discussions.


\section{Proof of Lemma 2}

In this section we prove Lemma 2.
Let $f:I\to I$ be a topologically mixing $C^2$ map
and $\mu\in\mathcal{M}_f$ hyperbolic with $h_{\mu} (f)>0$.
Without loss of generality we assume that $\mu$ is a linear combination
$\mu =c_1\nu_1 +\cdots + c_q\nu_q$ of 
$\nu_1,\ldots , \nu_q\in\mathcal{H}_f$
such that $c_i=l_i / l\enspace (i=1,\ldots , q)$ 
for some integers $l, l_1,\ldots , l_q \ge 1$ 
with $l_1+\cdots +l_q = l ,$
because any hyperbolic measure
is approximated by a measure of this form.
Let
$\varphi_1 ,\ldots , \varphi_p \in C(I)$
and $\varepsilon >0. $
For each $i=1,\ldots , q$ by Proposition 1
there are an integer $k_i \ge 1$,
a compact interval $L_i\subset I$ and
a family $\mathcal{K}_i$ of pairwise disjoint compact intervals
with $K\subset L_i =f^{k_i}(K) $ for all $K\in\mathcal{K}_i$ 
such that
$$(\log \sharp\mathcal{K}_i )/k_i\ge h_i  -\varepsilon /2 ,
\quad
| \frac{1}{k_i} \log\vert (f^{k_i})' (x) \vert 
-\lambda_i | \le\varepsilon /2$$
and  
$$|\frac{1}{k_i}S_{k_i} \varphi_j (x) - \int \varphi_j d\mu |\le
\varepsilon /2 \quad (j=1,\ldots , p)$$
for all $x\in\sqcup_{K\in\mathcal{K}_i } K$,  
where $h_i := h_{\nu_i}(f)$
and $\lambda_i := \int\log |f'|d\nu_i .$
We assume that
$k_1=k_2=\cdots =k_q$ and that
$f^{k_1} L_i =I$ holds for $i=1,\ldots ,q$
without loss of generality
by taking the iterations of $f$ if necessary. 
For each $i=1,\ldots , q$ and an integer $s\ge 1$ we set
$$\mathcal{P}_i(s) :=\Big\{ P=\bigcap_{j=0}^{s-1}f^{-jk_1}K_j :
 K_0, \ldots , K_{s-1}\in \mathcal{K}_i \Big\} .$$
Then for any $P\in \mathcal{P}_i(s)$
we have
$f^{sk_1} P=L_i,$ and hence
$f^{(s+1)k_1} P= I.$
Fix a large integer $n\ge 1$, and let
$$\mathcal{Q} :=\Big\{ Q=\bigcap_{i=1}^{q}f^{-m_i}P_i :
 P_i\in \mathcal{P}_i ( nl_i  ) \text{ for } i=1,\ldots , q  \Big\} ,$$
where  
$m_i: =  k_1 \sum_{j=1}^{i-1} (nl_j +1)$
for $i=1,\ldots , q.$
Remark that 
$f^{(nl+q)k_1} Q=I$ holds for each $Q\in\mathcal{Q} .$
Now we put
$$\varLambda :=\bigcap_{l=0}^{\infty} f^{-(nl+q)k_1}
 \Big( \bigsqcup_{Q\in\mathcal{Q} } Q \Big) .$$
Then it is a compact $f^{(nl+q)k_1}$-invariant set, and 
the restriction of $f^{(nl+q)k_1}$ to $\varLambda$ is 
semiconjugate to a full-shift of $\sharp \mathcal{Q}$-symbols.
Thus the topological entropy $h(f^{(nl+q)k_1}|_{\varLambda} )$ 
is estimated from below as follows:
\begin{align*}
h(f^{(nl+q)k_1}|_{\varLambda} ) 
&\ge \log \sharp \mathcal{Q} 
=\sum_{i=1}^q nl_i \cdot \log \sharp \mathcal{K}_i \\
&\ge \sum_{i=1}^q nl_i k_i(h_i -\varepsilon /2 )  
= nl k_1(h_{\mu }(f) - \varepsilon /2 ) , 
\end{align*}
and then, by the variational principle for entropy \cite{W}
we obtain
an $f$-invariant Borel probability measure $\nu$ with
$\nu (\cup_{i=0}^{(nl+q)k_1-1} f^i \varLambda )=1$ such that
\begin{align*}
h_{\nu} (f) 
&\ge \{h(f^{(nl+q)k_1}|_{\varLambda} ) - \varepsilon /4\} /
\{(nl+q)k_1\}  \\ 
&\ge nl k_1 (h_{\mu }(f) - 3\varepsilon /4 ) /\{(nl+q)k_1\} \\
&\ge  h_{\mu}(f) - \varepsilon 
\end{align*}
since $n\ge 1$ is large. 
The measure $\nu$ above can be taken as ergodic by
considering the ergodic decomposion if necessary.
Then by the Ruelle inequality  we have
$$\lambda_{\nu} (f)\ge h_{\nu} (f)\ge h_{\mu}(f) - \varepsilon >0,$$
and hence $\nu\in\mathcal{H}_f$, if $\varepsilon >0$ is small.
Also, 
\begin{align*}
\log |(f^{(nl+q)k_1})'(x)| 
&= \sum_{i=1}^q \log |(f^{(nl_i+1)k_1})'( f^{m_i} (x))  | \\
&\le \sum_{i=1}^q   \{ 
nl_i k_1 (\lambda_i + \varepsilon /2 )
+k_1\cdot \max_{y\in I}\log |f'(y) | \} \\
&=  nlk_1 (\int\log |f'|d\mu  + \varepsilon /2 ) 
+qk_1\cdot \max_{y\in I}\log |f'(y) | \\
&\le  (nl+q)k_1 (\int\log |f'|d\mu + \varepsilon  ) 
\end{align*}
holds whenever $x\in \varLambda$, because $n\ge 1$ is large.
Thus we have 
\begin{align*}
\lambda_{\nu} (f) &= \int\log |f'|d\nu \\
&\le \max_{x\in \varLambda}\frac{1}{(nl+q)k_1}
\log  |(f^{(nl+q)k_1})'(x) | \\
&\le \int\log |f'|d\mu +\varepsilon .
\end{align*}
Moreover, 
\begin{align*}
| S_{(nl+q)k_1} &\varphi_j (x) 
- (nl+q)k_1\int\varphi_j d\mu | \\
&\le \sum_{i=1}^q   
| S_{nl_ik_1} \varphi_j (f^{m_i}(x)) 
- nl_ik_1\int\varphi_j d\nu_i | 
+2qk_1 \max_{y\in I}|\varphi_j (y) |  \\
&\le \sum_{i=1}^q 
nl_i k_1\varepsilon /2 + 2qk_1 \max_{y\in I}|\varphi_j (y) |  \\
&\le (nl+q)k_1\varepsilon
\end{align*}
holds whenever $x\in \varLambda ,$ 
and from which
$$|  \int\varphi_j d\nu - \int\varphi_j d\mu | \le \varepsilon$$
for each $j=1,\ldots , p.$ 
This completes the proof of the lemma.

\section{Proofs of Theorems}

Let $f:I\to I$ be a 
topologically mixing $C^2$ map with 
a return time function $R: J\to \mbox{\boldmath $N$}\cup\{\infty \}$
satisfying the assumptions (H1)-(H4).
The proofs of Theorems 3 and 4 are reduced to 
two propositions mentioned below.
The one of them concerning with the variational principle
gives an upper estimate of the dimension spectra for Birkhoff averages.

\begin{prop}
Let $\varphi_1 ,\ldots , \varphi_p  \in C(I)$ and
$\alpha_1 , \ldots , \alpha_p \in \mbox{\boldmath $R$}$.
Then
\begin{align*}
\dim_H\Big\{ &x\in I :\lim_{n\to\infty} 
\frac{1}{n}S_n\varphi_j (x)  = \alpha_j   
\enspace (j=1,2,\ldots , p) \Big\} \\
\le \sup{}^+&\Big\{ D(\nu ) :\nu\in\mathcal{H}_f ,
| \int\varphi_j d\nu  -\alpha_j  | < \varepsilon 
\enspace (j=1,2,\ldots , p) \Big\} 
\end{align*}
holds whenever $\varepsilon >0.$

\end{prop}
The assumptions on the return time function is necessary 
for this proposition.
Another proposition obtained from Moran-like geometric constructions 
\cite{P2} gives
a lower estimate of the dimension spectra for Birkhoff averages.
\begin{prop}
Let $\varphi_1 ,\varphi_2 ,\ldots \in C(I)$
and $\alpha_1 ,\alpha_2, \ldots \in \mbox{\boldmath $R$} .$
Assume that a sequence $\nu_1 ,\nu_2 , \ldots \in\mathcal{H}_f$
satisfies $$\lim_{p\to\infty}\max_{1\le j\le p}
| \int\varphi_j d\nu_p  -\alpha_j  | = 0.$$ 
Then
$$\dim_H \Big\{ x\in I : \lim_{n\to\infty} 
\frac{1}{n}S_n\varphi_j (x) = \alpha_j \enspace (j=1,2,\ldots )  \Big\} \\
\ge \limsup_{p\to\infty} D( {\nu_p} ) .$$
\end{prop}

Theorem 3 follows immediately from Propositions 6 and 7
by considering
the sequences $\{ \varphi_j \}_{j=1}^{\infty} $
with $\varphi_j = \varphi \enspace (j\ge 1)$ and 
$\{ \alpha_j \}_{j=1}^{\infty} $ 
with
$\alpha_j =\alpha \enspace (j\ge 1)$ 
for given $\varphi\in C(I)$ and 
$\alpha \in \mbox{\boldmath $R$} .$
Theorem 4 is also obtained from the propositions as follows.
For $\mu\in\mathcal{M}_f$
put
$$D_0 : = \inf_{\mathcal{U}}\sup \{ D(\nu) : 
\nu\in\mathcal{H}_f\cap\mathcal{U} \} ,$$
where the infimum is taken over all of 
the neighborhoods $\mathcal{U}$ of $\mu$ in $\mathcal{M}$.
To show $\dim_H (B(\mu ))\ge D_0$ we assume that $D_0 >0$,
otherwise it is trivial.
Choose a sequence $\varphi_1 , \varphi_2 , \ldots \in C(I)$
to be dense in $C(I)$. 
Then 
for any $\gamma >0$ 
taking large $p\ge 1$ and small $\varepsilon >0$
we have
\begin{align*}
\dim_H &(B(\mu) ) \\
\le &\dim_H \Big\{ x\in I : 
\lim_{n\to\infty} 
 \frac{1}{n}S_n\varphi_j (x) = 
\int\varphi_j d\mu  \enspace (j=1,2,\ldots , p ) 
 \Big\} \\
\le &\sup \Big\{ D(\nu ) :\nu\in\mathcal{H}_f ,\enspace
| \int\varphi_j d\nu  -\int \varphi_j d\mu  |< \varepsilon 
\enspace (j=1,2,\ldots , p) \Big\} \\
\le &D_0 +\gamma 
\end{align*}
by Proposition 6.
Letting $\gamma \to 0$ we obtain 
\begin{equation}
\dim_H (B(\mu) ) \le D_0 .
\label{hdbu}
\end{equation} 
Moreover, taking a sequence $\nu_1 , \nu_2 ,\ldots\in\mathcal{H}_f$
such that
$$\lim_{p\to\infty} D(\nu_p) = D_0
\quad\text{ and }\quad
\lim_{p\to\infty}\max_{1\le j\le p}
| \int\varphi_j d\nu_p  -\int\varphi_j d\mu | = 0 $$  
we have
\begin{align}
\dim_H &(B(\mu) )   \label{hdbl} \\
=&\dim_H \Big\{ x\in I : \lim_{n\to\infty} 
\frac{1}{n}S_n\varphi_j (x) = \int\varphi_j d\mu \enspace (j=1,2,\ldots ) 
 \Big\}  \nonumber \\
\ge &\limsup_{p\to\infty} D( {\nu_p} ) \quad =D_0 \nonumber 
\end{align}
by Proposition 7.
Combining (\ref{hdbu}) and (\ref{hdbl})
we obtain 
$$\dim_H (B(\mu) ) = D_0  = \inf_{\mathcal{U}}\sup \{ D(\nu) : 
\nu\in\mathcal{H}_f\cap\mathcal{U} \} ,$$
and hence Theorem 4.
For the proof of Theorem 5 we need the following:
\begin{prop}
Let $\varphi\in C(I)$ and $\mu_1 ,\mu_2\in \mathcal{H}_f$
with $\displaystyle \int\varphi d\mu_1 \not= \int\varphi d\mu_2.$
Then 
$$\dim_H (I_{\varphi})  \ge \min\{ D(\mu_1) , D(\mu_2)   \} .$$
\end{prop}
Theorem 5 is obtained from Proposition 8 as follows. 
We claim that for any $\varepsilon >0$ there are 
$\mu_1 , \mu_2 \in \mathcal{H}_f$ with $\mu_1\not= \mu_2$
such that $D(\mu_i )\ge 1-\varepsilon$ for $i=1,2.$
In fact, 
taking $\mu_1\in\mathcal{H}_f$ 
as an absolutely continuous invariant probability measure
we have
\begin{align*}
D(\mu_1) =h_1 /\lambda_1  =1 ,
\end{align*}
where $h_1 := h_{\mu_1}(f)$ and 
$\lambda_1 :=\lambda_{\mu_1} (f)$, respectively.
To obtain $\mu_2$ fix a small number $\gamma >0$ such that 
$(h_1 -2\gamma )/ (\lambda_1 +2\gamma) \ge 1-\varepsilon$, and  
$\nu\in\mathcal{H}_f$ such that $\nu\not= \mu_1$.
Let
$\nu_t := t\nu +(1-t)\mu_1 \in\mathcal{M}_f$ for $0<t<1.$ 
Then $\nu_t$ is a hyperbolic measure 
with positive metric entropy 
such that $\nu_t\not= \mu_1$ for all $t$. 
Moreover, 
$$h_{\nu_t} (f) \ge h_1 - \gamma ,\quad 
\int \log |f'|d\nu_t \le \lambda_1  +\gamma$$ 
holds if $t>0$ is sufficiently small.  
By lemma 2 there is $\mu_2\in\mathcal{H}_f$ 
sufficiently close to $\nu_t ,$ 
and then
$\mu_2\not=\mu_1$,
such that 
$$h_2 \ge h_{\nu_t} (f) - \gamma \ge h_1 -2\gamma$$
and 
$$\lambda_2 \le \int \log |f'| d\nu_t +\gamma
\le \lambda_1  +2\gamma ,$$
where $h_2 := h_{\mu_2}(f)$ and 
$\lambda_2 :=\lambda_{\mu_2} (f)$, respectively.
Then
$$D(\mu_2 ) = h_2 /\lambda_2  
\ge (h_1-2\gamma)/(\lambda_1+2\gamma ) \ge 1-\varepsilon ,$$ 
and thus, the claim is established.
Taking a continuous function $\varphi :I\to \mbox{\boldmath $R$}$
such that 
$\displaystyle \int\varphi d\mu_1 \not= \int\varphi d\mu_2 ,$
we have
\begin{align*}
\dim_H (I\setminus QR(f) ) 
&\ge \dim_H (I_{\varphi})  \\
&\ge \min\{ D(\mu_1) , D(\mu_2)   \} \\
&\ge 1-\varepsilon 
\end{align*}
by Proposition 8.
Letting $\varepsilon \to 0$ we obtain Theorem 5.

\vspace{0.2cm}

\noindent {\bf Remark.}
The assertion of Theorem 5 is valid without assuming 
the existence of a return time function
satisfying (H1)-(H4), if the map $f: I\to I$ has an  
absolutely continuous invariant probability measure.

\section{Proof of Proposition 6}

We prepare a lemma for the proof of Proposition 6. 
Let $(Y, \mathcal{B} , m)$ be a finite measure space and
$Y_1 , \ldots , Y_l \in\mathcal{B} $ 
pairwise disjoint subsets of $Y$ with 
positive measures.
We consider a measurable map $g: \sqcup_{j=1}^l Y_j \to Y$ 
satisfying the following properties:
\begin{enumerate}
\item for each $j=1,\ldots , l, $
$g_j :=g|_{Y_j} : Y_j\to Y$ is a bi-nonsingular bijection;
\item for any sequence $\{  a_j \}_{j=0}^{\infty}$ with 
$a_j\in \{ 1,\ldots ,l \} \enspace (j=0,1,2,\ldots )$ 
the set $\cap_{n=0}^{\infty} Y_{a_0\ldots a_{n-1}}$ 
consists of a single point,
where 
$$Y_{a_0\cdots a_{n-1}} :
= Y_{a_0}\cap g^{-1}Y_{a_1}\cap\cdots\cap g^{-(n-1)}Y_{a_{n-1}};$$
\item there are a version ${\rm Jac} (g) >0$ 
of the Radon-Nikodym derivative $\frac{dm\circ g}{dm}$ 
such that
$$ \lim_{n\to\infty}\sup_{J\in W_n}\sup_{x, y\in Y_J}
| \log {\rm Jac}(g)(x)-\log {\rm Jac}(g)(y) | =0$$
and $C\ge 1$, a distortion constant, such that
for any integer $n\ge 1$ and
$J \in W_n$
$$\frac{\prod_{i=0}^{n-1} {\rm Jac} (g) (g^i(x)) }
{\prod_{i=0}^{n-1} {\rm Jac} (g) (g^i(y)) }
\le C$$
holds whenever $x, y\in Y_J,$
where 
$W_n$ stands for the set of the words 
$J=(a_0 \cdots a_{n-1})$
of length $n$ with 
$a_i\in \{ 1,2,\ldots , l \}$ for each $i=0,1,\ldots ,n-1.$
\end{enumerate}
Then we call $g:\varLambda \to\varLambda$ a 
{\it finite Markov system} induced by $( Y, \{ Y_i\}_{i=1}^l , g) ,$
where $\varLambda := \cap_{n=0}^{\infty} g^{-n} (\sqcup_{j=1}^l Y_j) .$
It is isomorphic to a full-shift of $l$-symbols, 
and the space of the probability measures supported on $\varLambda$
is compact.
The following lemma obtained from a standard argument 
on the variatinal principle for pressure \cite{W}
is a generalization of Lemma 7 of \cite{C4}.

\begin{lem}
For any $\delta \ge 0$ 
a finite Markov system $g:\varLambda \to \varLambda$ 
induced by $(Y, \{ Y_j\}_{j=1}^l , g)$
with a distortion constant $C\ge 1$ has an
invariant ergodic probability measure $\mu$ on $\varLambda$ such that
$$h_{\mu} (g) -\delta \int\log {\rm Jac} (g) d\mu \ge
\log \sum_{j=1}^l m(Y_j)^{\delta} - \delta \log (C m(Y) ) .$$
\end{lem}

\noindent
{\it Proof.}
For each $J=(a_0\ldots a_{n-1})\in W_n$  and 
$ p= 1,\ldots , l$
we write $Y_{Jp} := Y_J\cap g^{-n}Y_p .$ 
Then since
\begin{align*}
\frac{m(Y_{Jp})}{m(Y_{J})} 
\ge \frac{\prod_{i=0}^{n-1}\inf_{x \in Y_{a_0\ldots a_{n-1}p}}
{\rm Jac} (g)(g^i(x))^{-1} m(Y_p)}
{\prod_{i=0}^{n-1}\sup_{y\in Y_{a_0\ldots a_{n-1}}}
{\rm Jac} (g)(g^i(y))^{-1} m(Y)} 
\ge C^{-1} \frac{m(Y_p)}{m(Y)} 
\end{align*}
for each $p=1,\ldots , l,$
we have 
$$\frac{\sum_{p=1}^{l} m(Y_{Jp})^{\delta}} {m(Y_{J})^{\delta}} \ge 
 \frac{\sum_{p=1}^{l} m(Y_p)^{\delta}}{(Cm(Y))^{\delta}}  
,$$
and hence
\begin{align*}
\sum_{J\in W_n}  m(Y_J)^{\delta}  
\ge &\frac{\sum_{p=1}^l m(Y_p)^{\delta} }{(Cm(Y))^{\delta}} \cdot
\sum_{J\in W_{n-1}}  m(Y_{J})^{\delta} \\  
\cdots \ge &\left\{
\frac{ \sum_{p=1}^l m(Y_p)^{\delta} }{(Cm(Y))^{\delta}} \right\}^{n-1}
\cdot\sum_{J\in W_1} m(Y_J)^{\delta} \\
= &\left\{ 
\frac{ \sum_{p=1}^l m(Y_p)^{\delta}}{(Cm(Y))^{\delta}} \right\}^{n}
\cdot (Cm(Y))^{\delta} .
\end{align*}
Thus we obtain
\begin{align}
\liminf_{n\to\infty}
\frac{1}{n} \log\sum_{J\in W_n} m(Y_J)^{\delta}
\ge \log \sum_{j=1}^l m(Y_j)^{\delta}  -\delta \log (Cm(Y) )
 .\label{eq:fms1}
\end{align}
Also, taking $x_J\in Y_J\cap\varLambda$ 
for each $J\in W_n $ 
we obtain a sequence of probability measures supported on $\varLambda$
by
$\displaystyle \mu_n :=\frac{1}{Z_n} \sum_{J\in W_n} m(Y_J)^{\delta}
 \delta_{x_J}^n $
for each $n\ge 1$,
where $Z_n := \sum_{J\in W_n} m(Y_J)^{\delta}$ and 
$\delta_{x_J}^n := (\delta_{x_J} +\delta_{g(x_J)} + 
\cdots + \delta_{g^{n-1}(x_J)} )/n .$
Then an accumulation point  
$\mu$ of the sequence $\{ {\mu}_n \}_{n=1}^{\infty}$
is a $g$-invariant probability measure 
supported on $\varLambda .$
Moreover,
\begin{align*}
\log &\sum_{J\in W_n} m(Y_J)^{\delta} 
=\log Z_n \\
=
&\sum_{J\in W_n} \mu_n (Y_J) 
\left\{ - \log \mu_n (Y_J) + \delta\log m(Y_J)   \right\} \\
\le
&\sum_{J\in W_n} \mu_n (Y_J) 
\left\{  - \log \mu_n (Y_J) + 
\delta\log \left( C m(Y)\prod_{i=0}^{n-1}
\text{Jac} (g)(g^i(x_J))^{-1} \right)  \right\} \\
=&-\sum_{J\in W_n} \mu_n (Y_J)\log \mu_n (Y_J) 
- n\delta\int\log {\rm Jac} (g) d\mu_{n} +\delta\log ( C m(Y) ) ,
\end{align*}
and hence
\begin{equation}
\limsup_{n\to\infty}\frac{1}{n} \log\sum_{J\in W_n} m(Y_J)^{\delta}
\le h_{\mu} (g) - \delta\int\log {\rm Jac} (g) d\mu . \label{eq:fms2}
\end{equation}
Combining (\ref{eq:fms1}) and (\ref{eq:fms2})
we obtain
$$h_{\mu} (g) -\delta\int\log {\rm Jac} (g) d\mu \ge
\log \sum_{j=1}^l  m(Y_j)^{\delta} - \delta \log (C m(Y)).$$
The above $\mu$ can be taken as an ergodic measure
by considering the ergodic decomposition if necessary.
We have proved the lemma. 

\vspace{0.2cm}

Now we prove Proposition 6.
For given $\varphi_1 ,\ldots , \varphi_p  \in C(I),$
$\alpha_1 , \ldots , \alpha_p \in \mbox{\boldmath $R$}$
and $\varepsilon >0$ we set
$$\delta_0 := \sup{}^+ \Big\{ D(\nu) : \nu\in\mathcal{H}_f ,
|\int\varphi_j d\nu -\alpha_j  | < \varepsilon 
\enspace (j=1,2,\ldots , p) \Big\} \in [0,1] .$$
and show that
\begin{equation}
\label{ineq:delta}
\dim_H\Big\{ x\in I :\lim_{n\to\infty} 
\frac{1}{n}S_n\varphi_j (x) = \alpha_j  
\enspace (j=1,2,\ldots , p) \Big\} 
\le \delta_0 .
\end{equation}
Remark that 
\begin{equation}\label{h-dl}
h_{\nu} (f) -\delta_0 \lambda_{\nu} (f) \le 0
\end{equation}
holds whenever
$\nu\in\mathcal{H}_f$ satisfies 
$\displaystyle | \int\varphi_j d\nu  -\alpha_j | < \varepsilon \enspace 
(j=1,2,\ldots , p)$
by the definition of $D (\nu )$.
Since $f: I\to I$ is topologically mixing, 
taking an integer $k\ge 1$ such that $f^k(J)=I$ 
for the base set $J$ of the return time function,
we have
\begin{align*}
f^k \Big(\Big\{ &x\in J : 
\lim_{n\to\infty}  \frac{1}{n}S_n\varphi_j (x) = \alpha_j 
\enspace (j=1,2,\ldots , p) \Big\} \Big) \\
&=  \Big\{ x\in I : 
\lim_{n\to\infty}  \frac{1}{n}S_n\varphi_j (x) = \alpha_j 
\enspace (j=1,2,\ldots , p) \Big\} ,
\end{align*}
and from which it follows that
\begin{align*}
\dim_H \Big\{ &x\in I : 
\lim_{n\to\infty}  \frac{1}{n}S_n\varphi_j (x) = \alpha_j 
\enspace (j=1,2,\ldots , p) \Big\} \\
=  &\dim_H \Big\{ x\in J : 
\lim_{n\to\infty}  \frac{1}{n}S_n\varphi_j (x) = \alpha_j 
\enspace (j=1,2,\ldots , p) \Big\}  ,
\end{align*}
because $f^k :I\to I$ is smooth.
Moreover,
$$\dim_H \Big\{ x\in J :  
\lim_{n\to\infty}  \frac{1}{n} S_n\varphi_j (x) =\alpha_j 
\enspace (j=1,2,\ldots , p) \Big\}
\subset \bigcup_{N=1}^{\infty }\varGamma_N $$
holds,
where $$\varGamma_N :=\Big\{ 
x\in J :  
| \frac{1}{n}S_n\varphi_j (x) -\alpha_j |< \varepsilon /2
\enspace (j=1,2,\ldots , p , \enspace n\ge N ) \Big\} .$$
Thus for the inequality (\ref{ineq:delta})
it is enough to show that 
\begin{equation} \label{eq:d0}
H^{\delta_0} (\varGamma_N  ) <\infty
\end{equation}
for each integer $N\ge 1$,
where $H^{\delta_0}$ denotes the 
$\delta_0$-dimensinal Hausdorff outer measure.
For an integer $n\ge 1$ let $\mathcal{A}_n$ 
be a finite partition of the base set $J$
such that $A\in\mathcal{A}_n$ iff
$A$ is a connected component of
\begin{align*}
\{ x \in J :&R(x)= k_1 ,\enspace R(f^{k_1}(x))=k_2, \enspace 
\ldots \enspace ,\\
&R(f^{k_1+\cdots +k_{l-2}}(x))= k_{l-1},\enspace   
R(f^{k_1+\cdots +k_{l-1}}(x))\ge k_l  
\}
\end{align*} 
for some integers $k_1 , \ldots ,  k_l\ge 1$  
with $k_1+\cdots + k_l=n$,
and $\mathcal{D}_n$ a  family of
compact intervals such that 
$B\in\mathcal{D}_n$ iff
$B$ is a connected component of
$$\{ x \in J :R(x)= k_1 ,\enspace R(f^{k_1}(x))=k_2,
 \enspace \ldots \enspace ,
\enspace R(f^{k_1+\cdots +k_{l-1}}(x))= k_l \} $$
for some integers $k_1 , \ldots ,  k_l\ge 1$  
with $k_1+\cdots + k_l=n .$  
Set 
\begin{align*}
\mathcal{B}_n := \Big\{ A\in\mathcal{A}_n : 
| \frac{1}{n}S_n\varphi_j (x_A)-\alpha_j |< \varepsilon /2 
\enspace &(j=1,2,\ldots , p  ) \\
&\text{ for some } x_A\in A \Big\} .
\end{align*}
Then we have 
\begin{equation} \label{g-a} 
\varGamma_N\subset 
\bigcap_{n =N}^{\infty} \bigsqcup_{A\in\mathcal{B}_n} A .
\end{equation}
Putting 
$$\mathcal{B}^*_{n, s} := 
\{  B\in\mathcal{D}_{n+s} : B\subset A \text{ for some }
A\in\mathcal{B}_n \}$$
for each $s\ge 1$,
we obtain
$$|A|\le  C{\gamma_0}^{-1}
\sum_{s=0}^{l_0-1} \sum_{B\in \mathcal{B}^*_{n, s} , B\subset A} |B|
$$
for each $A\in\mathcal{B}_n $
from the assumptions (H3) and (H4) for the return time function.
Then, since $0\le \delta_0 \le 1$ we have
\begin{align*}
\sum_{s=0}^{l_0-1} \sum_{B\in \mathcal{B}^*_{n, s}}
|B|^{\delta_0} 
&= 
\sum_{A\in\mathcal{B}_n}\sum_{s=0}^{l_0-1} 
\sum_{B\in \mathcal{B}^*_{n, s}, B\subset A} 
|B|^{\delta_0}  \\
&\ge 
\sum_{A\in\mathcal{B}_n} \Big(\sum_{s=0}^{l_0-1} 
\sum_{B\in \mathcal{B}^*_{n, s}, B\subset A} 
|B|\Big)^{\delta_0}  \\
&\ge
C^{-\delta_0}{\gamma_0}^{\delta} 
\sum_{A\in\mathcal{B}_n}
|A|^{\delta_0}  ,
\end{align*}
and hence
\begin{equation}
\label{ld1}
\sum_{B\in \mathcal{B}^*_{n, s}}
|B|^{\delta_0} 
\ge
{l_0}^{-1}C^{-\delta_0}{\gamma_0}^{\delta_0} 
\sum_{A\in\mathcal{B}_n}
|A|^{\delta_0}  
\end{equation}
holds some integer $s=s(n)$ with $0\le s\le l_0-1 .$
Also, from the definition of $\mathcal{B}^*_{n, s}$
we have
\begin{equation}
\label{ineq:n+s}
|\frac{1}{n+s}S_{n+s}\varphi_j (x)-\alpha_j | < \varepsilon 
\enspace (j=1,2,\ldots , p  ) 
\end{equation}
for all  $x\in B$ with $B\in \mathcal{B}^*_{n, s}$
if $n$ is large.
Since 
$(  J, \mathcal{B}^*_{n, s} , f^{n+s}  )$
is a finite Markov system with the distortion constant $C$,
by Lemma 9
there is an $f^{n+s}$-invariant ergodic probability measure
$\mu_n$ supported on 
$\varLambda_n :
=\cap_{k=0}^{\infty} f^{-(n+s)k}(\sqcup_{B\in \mathcal{B}^*_{n, s} } B )$
such that
\begin{equation*} 
h_{\mu_n} (f^{n+s} ) -\delta_0 \int \log |(f^{n+s})'| d\mu_n
\ge \log \sum_{B\in \mathcal{B}^*_{n, s}} |B|^{\delta_0} 
- \delta_0 \log(C |J|) .
\end{equation*}
Let 
$${\nu}_n := \frac{1}{n+s}\sum_{i=0}^{n+s-1}\mu_n\circ f^{-i} 
\in\mathcal{M}_f .$$
Then it is ergodic, 
and supported on  
$X_n :=\cup_{i=0}^{n+s-1} f^i (\varLambda_n ) .$
The set $X_n$ is a repeller for $f$
by the assumption (H1), and hence $\nu_n\in \mathcal{H}_f .$
Moreover, by (\ref{ineq:n+s}) we have
\begin{equation*}
| \int\varphi_j d\nu_n -\alpha_j |
\le \max_{B\in \mathcal{B}^*_{n, s}}\max_{x\in B}
| \frac{1}{n+s} S_{n+s}\varphi_j (x) -\alpha_j | 
< \varepsilon 
\end{equation*}
for all $j=1,2,\ldots , p ,$
and then by (\ref{h-dl}) we obtain
\begin{align}
\label{ld2}
\log \sum_{B\in\mathcal{B}^*_{n, s}} |B|^{\delta_0}
&\le 
 h_{{\mu}_n} (f^{n+s} ) -\delta_0 \int \log |(f^{n+s})'| d {\mu}_n
+\delta_0\log (C |J|) \\
&=
(n+s) ( h_{\nu_n} (f) -\delta_0 \int \log |f'| d\nu_n)
+\delta_0\log (C |J|)  \nonumber \\
&\le \delta_0\log (C |J|).  \nonumber 
\end{align}
Since $\max_{A\in  \mathcal{B}_n} |A| \le \varepsilon_n \to 0$
as $n\to\infty$ by the assumption (H2), 
combining (\ref{g-a}), (\ref{ld1}) and  (\ref{ld2}) we have
\begin{align*}
H^{\delta_0} (\varGamma_N )
&\le \limsup_{n\to\infty} \sum_{A\in \mathcal{B}_n}| A|^{\delta_0} \\
&\le  l_0 ( C^2 {\gamma_0}^{-1} |J|)^{\delta_0} 
\quad <\infty
\end{align*}
holds for all $N\ge 1 .$
We have thus obtained (\ref{eq:d0}), 
and hence Proposition 6.

\section{Proof of Proposition 7}

Let $\{\varphi_j \}_{j=1}^{\infty}\subset C(I), 
\{\alpha_j \}_{j=1}^{\infty}\subset\mbox{\boldmath $R$} $
and  $\{\nu_j \}_{j=1}^{\infty}\subset \mathcal{H}_f$ 
be as in Proposition 7.
Taking a subsequences if necessary, we may assume that 
$D_p:= D(\nu_p)\enspace (p=1,2,\ldots )$ is a sequence of 
positive numbers which converges to some number $D_0\in [0,1]$ 
as $p\to \infty$ and 
$$\beta_p :=
\max_{1\le j\le p}
| \int\varphi_j d\nu_p  -\alpha_j  | \in (0, 1) 
\enspace (p=1,2,\ldots )$$
a decreasing sequence to zero,
respectively.
Then for each $p\ge 1$ we have  $\lambda_p \ge h_p >0$ holds
by the Ruelle inequality,
where
$ \lambda_p := \lambda_{\nu_p} (f)$
and $ h_p := h_{\nu_p}(f) .$
Take a sequence $\{ \gamma_p \}_{p=1}^{\infty}$ of
positive small numbers such that
$$\lim_{p\to\infty} \frac{h_p-\gamma_p}{\lambda_p+2\gamma_p} =
D_0 .$$
Then for each integer $p\ge 1$ by Proposition 1
there are an integer $k_p \ge 1$,
a compact interval $L_p\subset I$ and
a family $\mathcal{K}_p$ of pairwise disjoint compact intervals
with $K\subset L_p =f^{k_p}(K) $ for all $K\in\mathcal{K}_p$ 
such that
$$(\log \sharp\mathcal{K}_p )/k_p\ge h_p -\gamma_p   ,
\quad
| \frac{1}{k_p} \log\vert (f^{k_p})' (x) \vert 
-\lambda_p | \le  \gamma_p  $$
and  
$$|\frac{1}{k_p}S_{k_p} \varphi_j (x) - \int \varphi_j d\mu_p |\le
\beta_p  \quad (j=1,\ldots , p)$$
for all $x\in\sqcup_{K\in\mathcal{K}_p } K .$  
Taking some iterations of the map $f$ if neccesary,
we assume that $f^{k_p} (L_p)= I$
without loss of generality.
Choose a sequence $\{ q_p \}_{p=1}^{\infty }$ of 
positive integers such that $q_p$ is so large
as compared with $q_1,\ldots , q_{p-1}, k_1,\ldots , k_{p+1}$ 
for each $p\ge 1$,
and put 
$n_p :=(q_1+1)k_1+\cdots +(q_p+1)k_p$ 
for $p\ge 0 .$
For an integer $l\ge 1$ take a pair of integers $p=p(l)\ge 0$ 
and $s=s(l)$ with $1\le s\le q_{p+1}$ such that 
$l=q_1+\cdots + q_{p} + s$, and set
\begin{align*}
\mathcal{S}(l) :=\{ (K_1,\ldots , K_l) :
K_i\in \mathcal{K}_r &\text{ if } t_{r-1}< i\le t_r \\
&\text{ for } i=1,\ldots ,l \}
\end{align*}
where $t_r=q_1+\cdots + q_r$ for each integer $r\ge 0 .$
For each $(K_1,\ldots , K_l)\in \mathcal{S}(l) $
take an interval $Q(K_1,\ldots , K_l)$ inductively as follows.
Let $Q(K) =K$ for $K\in\mathcal{K}_1$,  
and choose a compact interval 
$Q(K_1,\ldots , K_{l}) \subset  Q(K_1,\ldots , K_{l-1})$
such that
 $$f^{n_{p-1}+ (s-1)k_p } Q(K_1,\ldots , K_l) =K_l$$ 
if $Q(K_1,\ldots , K_{l-1})$ has been defined.
Set 
$$\mathcal{Q} (l) := \{ Q:=Q(K_1,\ldots , K_l) :
(K_1,\ldots , K_l) \in \mathcal{S}(l)  \} .$$
Then it holds that
\begin{align} \label{sharpQ}
\sharp \mathcal{Q}(l)
&=\sharp \mathcal{S}(l)  \quad
= (\sharp\mathcal{K}_1)^{q_1}\cdots (\sharp\mathcal{K}_{p})^{q_{p}}
\cdot (\sharp\mathcal{K}_{p+1})^s \\
&\ge  (\sharp\mathcal{K}_{p})^{q_{p}}
\cdot (\sharp\mathcal{K}_{p+1})^s
 \nonumber \\
&\ge \exp \{ q_{p}k_{p}(h_{p}-\gamma_{p}  ) 
+  sk_{p+1}(h_{p+1}-\gamma_{p+1}  ) \} \nonumber .
\end{align}
Put $m_l := n_{p} + (s+1)k_{p+1}$, where integers $p\ge 1$ and 
$s$ with $1\le s\le q_p$ are chosen so that
$q_1+\cdots + q_{p-1}+s =l$.
Then for each $Q\in\mathcal{Q} (l),$
$$f^{m_l}(Q)=f^{k_{p+1}}(L_{p+1}) =I$$ 
holds by the definition of $\mathcal{Q} (l) .$
And,
since $q_p$ is so large as compared with 
$q_1,\ldots , q_{p-1}$ and $k_1, \ldots , k_{p+1}$ we have
\begin{align*}
|(&f^{m_l} )' (x) |
= \Big( \prod_{i=1}^{p} |(f^{(q_i+1)k_i})'(f^{n_{i-1}}(x)) |\Big)
\cdot |(f^{(s+1)k_{p+1}})'(f^{n_{p}}(x)) | \\
&\le  \exp\Big\{
\sum_{i=1}^{p}q_ik_i(\lambda_i  +\gamma_i   )  
+ sk_{p+1}( \lambda_{p+1} +\gamma_{p+1} ) \Big\} \cdot 
\max_{y\in I}|f'(y) |^{k_1+\cdots +k_{p+1}} \\
&\le  \exp\{
q_{p}k_{p}(\lambda_{p} + 2\gamma_{p} )   
+ sk_{p+1} ( \lambda_{p+1} + 2\gamma_{p+1} )  \} 
\end{align*}
holds for all $x\in Q.$
Then by the mean value theorem we get
\begin{equation}
\label{mvt}
| Q| \ge |I|\cdot
\exp \{-q_{p}k_{p}(\lambda_{p}+ 2\gamma_{p}  )
-sk_{p+1}(\lambda_{p+1}+ 2\gamma_{p+1}  )       \}
\end{equation}
for all $Q\in\mathcal{Q}(l)$
where $p\ge 0$ and $s$ with $1\le s\le q_{p+1}$ 
are integers such that $q_1 +\cdots +q_{p} + s =l .$
Moreover, for $j=1,2,\ldots , p$, we obtain
\begin{align*}
|S_{m_l}\varphi_j (x) - &m_l\alpha_j | \\
\le
\sum_{i=1}^{p} \Big\{
|&S_{q_ik_i}\varphi_j (f^{n_{i-1}}(x)) - q_ik_i\int\varphi_j d\nu_i | \\
&+ | q_ik_i\int\varphi_j d\nu_i - q_ik_i\alpha_j | 
+ 2k_i \max_{y\in I} |\varphi_j (y)-\alpha_j  | \Big\} \\
+ |&S_{sk_{p+1}}\varphi_j (f^{n_{p}}(x)) 
- sk_{p+1}\int\varphi_j d\nu_{p+1} | \\
&+ | sk_{p+1} \int\varphi_j d\nu_{p+1} - sk_{p+1}\alpha_j | 
+ 2k_{p+1} \max_{y\in I} |\varphi_j (y)-\alpha_j  | \\
\le 2 \Big(
 \sum_{i=1}^{p} &q_ik_i \beta_i +  sk_{p+1}\beta_{p+1}
+ \max_{y\in I} |\varphi_j (y)-\alpha_j  |\sum_{i=1}^{p+1} k_i \Big) \\
\le 4 ( q_pk_p &\beta_{p}  + sk_{p+1}\beta_{p+1} ) \quad
\le 4 m_l \beta_{p}  ,
\end{align*}
and hence
\begin{equation}
\label{ineq:ep-1}
|\frac{1}{m_l} S_{m_l}\varphi_j (x) - \alpha_j | \le 4 \beta_{p}
\end{equation}
holds for all $x\in Q$ whenever $Q\in\mathcal{Q} (l) .$
Let
$$\varGamma :=\bigcap_{l=1}^{\infty}\bigsqcup_{Q\in\mathcal{Q}(l)} Q .$$
Then by (\ref{ineq:ep-1}) for any integer $j\ge 1$ 
$$\lim_{n\to\infty}\frac{1}{n}S_n\varphi_j (x)=\alpha_j$$
holds whenever $x\in \varGamma .$
Thus, it is enough to show 
\begin{equation*} 
\dim_H (\varGamma )\ge D_0
\end{equation*}
for the proof of the proposition,
and it follows from the existence of a probability measure
$\mu_0$  with $\mu_0 (\varGamma )>0$ such that
\begin{equation} 
\label{fro}
\liminf_{r\to 0} \frac{\log \mu_0 ([x-r, x+r])}{\log r} 
\ge D_0
\end{equation}
for all $x\in\varGamma .$
For each integer $l\ge 1$ and $Q\in\mathcal{Q}(l)$ take a point 
$x_Q\in Q\cap\varGamma$ and define a probability measure
$\mu_l :=\sum_{Q\in \mathcal{Q}(l)} \delta_{x_Q} /\sharp\mathcal{Q}(l) $ 
on $I$.
Then, $\mu_l (\varGamma ) =1$ holds.
Let $\mu_0$ be an accumulation point of $\{\mu_l \}_{l=1}^{\infty}$ 
in $\mathcal{M}$. 
Notice that $\mu_0 (\varGamma ) =1$ and that
for any $Q\in \mathcal{Q} (l)$,
\begin{align*}
\mu_n (Q)
= \lim_{n\to\infty}  \sharp \{ R \in \mathcal{Q}(n) : R\subset Q \} /
\sharp \mathcal{Q}(n) 
= 1/ \sharp \mathcal{Q}(l),
\end{align*}
holds whenever $n\ge l$, 
and hence $$\mu_0 (Q) = 1/ \sharp \mathcal{Q}(l) .$$
Now we show that the inequality (\ref{fro}) holds for $x\in\varGamma .$
For any small $r>0$
take a pair of integers $p\ge 0$ and $s$ with $1\le s\le q_{p+1}$ 
such that
\begin{align*}
|I|\cdot\exp \{-q_{p}k_{p}&(\lambda_{p}+ 2\gamma_{p}  )
-(s+1)k_{p+1}(\lambda_{p+1}+ 2\gamma_{p+1}  ) \}
< r \\
&\le |I|\cdot
\exp \{-q_{p}k_{p}(\lambda_{p}+ 2\gamma_{p}  )
-sk_{p+1}(\lambda_{p+1}+ 2\gamma_{p+1}  ) \} ,
\end{align*}
and let $l = q_1+\cdots + q_{p}+ s. $
Then by (\ref{mvt})
for any interval $Q\in \mathcal Q(l)$  
$$|Q| \ge |I|\cdot\exp \{-q_{p}k_{p}(\lambda_{p}+ 2\gamma_p  )
-sk_{p+1}(\lambda_{p+1}+ 2\gamma_{p+1}  ) \} \ge r ,$$
and hence the interval $[x-r , x+r ]$ intersects 
with at most 3 elements of $Q(l)$.
This implies that
\begin{align*}
\mu_0 ([x-r, x+r] )
&\le 3/\sharp\mathcal{Q}(l) \\
&\le 3  \exp \{ -q_{p}k_{p}(h_{p}-\gamma_{p}  ) 
-  sk_{p+1} (h_{p+1}-\gamma_{p+1} ) \} 
\end{align*}
by (\ref{sharpQ}),
and from which
\begin{align*}
\frac{\log \mu_0 ([x-r, x+r])}{\log r}  
&\ge 
\frac{ q_{p}k_{p}(h_{p}-\gamma_{p})
+ sk_{p+1} (h_{p+1}-\gamma_{p+1}) -\log 3}
{ q_{p}k_{p}(\lambda_{p}+ 2\gamma_{p}  )
+(s+1)k_{p+1}(\lambda_{p+1}+ 2\gamma_{p+1}  )-\log |I| } .
\end{align*}
Since $p\to \infty $ as  $r\to 0$ 
we have 
$$
\liminf_{r\to 0} \frac{\log \mu_0 ([x-r, x+r])}{\log r}  
\ge \lim_{p\to\infty} 
\frac{h_p -\gamma_p}{\lambda_p + 2\gamma_p}
= D_0 .$$
We have thus obtained the inequality (\ref{fro}), 
and hence Proposition 7.

\section{Proof of Proposition 8}

Let $\varphi\in C(I), \mu_1, \mu_2\in\mathcal{H}_f$ and assume that
$$\beta := | \int\varphi d\mu_1 -\int\varphi d\mu_2 | >0 .$$
Then by Proposition 1 
for any $\varepsilon >0$ and $i=1, 2$ 
there are an integer $k_i\ge 1$,
a compact interval $L_i\subset I$ and
a family $\mathcal{K}_i$ of pairwise disjoint compact intervals
with $K\subset L_i =f^{k_i}(K)$ for each $K\in\mathcal{K}$ 
such that
$$(\log \sharp\mathcal{K}_i )/k_i \ge h_i  - \varepsilon ,
\quad 
 |\frac{1}{k_i} \log | (f^{k_i} )' (x) | - \lambda_i |\le \varepsilon$$
and  
$$|\frac{1}{k_i}S_{k_i} \varphi (x) - \int \varphi d\mu_i |\le
\beta /8$$
for all $x\in\sqcup_{K\in\mathcal{K}_i } K,$
where $h_i := h_{\mu_i} (f)$ and 
$\lambda_i := \lambda_{\mu_i }(f)$ for $i=1,2,$
respectively.
Taking some iterations of the map $f$ if necessary,
we assume that $k_1 =k_2$ and 
that $f^{k_1}L_i=I$ for $i=1, 2$ without loss of generality.
Choose an increasing sequence $\{ q_p \}_{p=1}^{\infty} $ 
of positive integers 
such that $q_p$ is so large as compared with $q_1, \ldots , q_{p-1}$
for $p\ge 2$,
and set $n_{p} :=k_1 \sum_{i=1}^{p} (q_i+1)$ for each integer $p\ge 0$.
For an integer $l\ge 1$ there is a pair of integers 
$p =p(l)\ge 0$ and $s=s(l)$ with $1\le s\le q_{p+1}$ such that
$l=q_1+\cdots + q_{p}+s$. Then we put
\begin{align*}
\mathcal{S}(l) :=\{
(K_1,\ldots , K_l) :
K_i\in \mathcal{K}_1 &\text{ if } 
t_{2m}  < i \le  t_{2m+1} ,\\
K_i\in \mathcal{K}_2 &\text{ if } 
t_{2m+1} < i \le  t_{2(m+1)}  \\
&\text{ for some } m=0,1,2,\ldots \} ,
\end{align*}
where $t_n :=\sum_{j=1}^{n}q_j$ for $n\ge 0$,
and take a family 
$$\mathcal{Q}(l):= 
\{ Q(K_1,\ldots , K_l) :(K_1,\ldots , K_l) \in \mathcal{S}(l) \}$$
of pairwise disjoint compact intervals
$Q(K_1, \ldots , K_l)$
inductively chosen 
as in the proof of Proposition 7 so that
$Q(K_1,\ldots , K_{l})\subset Q(K_1,\ldots , K_{l-1})$
and
$f^{n_{p}+(s-1)k_1} (Q(K_1,\ldots , K_{l}))= K_l$
hold for each $(K_1,\ldots , K_{l})\in\mathcal{S}(l). $ 
Then it holds that
\begin{align*}
\sharp\mathcal{Q} (l) 
&= \sharp\mathcal{S} (l) \quad =(\sharp\mathcal{K}_1 )^{q_1}\cdots 
(\sharp\mathcal{K}_{p} )^{q_{p}}\cdot (\sharp\mathcal{K}_{p+1} )^{s} \\
&\ge (\sharp\mathcal{K}_{p} )^{q_{p}}\cdot (\sharp\mathcal{K}_{p+1} )^{s} \\
&\ge \exp [k\{q_{p} (h_{p} -\varepsilon ) + s (h_{p+1} -\varepsilon ) \}] ,
\end{align*}
where $\mathcal{K}_{2m-1} := \mathcal{K}_{1}, 
\mathcal{K}_{2m} := \mathcal{K}_{2} ,
h_{2m-1}:= h_1$ and $h_{2m}:=h_2$, respectively for each 
integer $m\ge 1.$
For any $Q\in \mathcal{Q}(l)$
we have $f^{n_{p}+(s+1)k_1}(Q)= I.$ 
And 
since $q_{p}$ is so large as compared with $q_1,\ldots , q_{p-1}$
we get
\begin{align*}
| (f^{n_{p}+(s+1)k_1} &)' (x) | 
= \prod_{i=1}^{p}  | (f^{(q_i+1)k_1 })' (f^{n_{i-1}}(x)) | 
\cdot
| (f^{(s+1)k_1})' (f^{n_{p}}(x)) |  \\
\le \exp&\Big[ k_1\Big\{ \sum_{i=1}^{p} q_i (\lambda_i +\varepsilon )
+s (\lambda_{p+1} +\varepsilon ) 
+ (p+1) \max_{y\in I}\log |f'(y) | \Big\}\Big] \\
\le \exp &[ k_1 \{ q_{p} (\lambda_{p} +2 \varepsilon )
+s (\lambda_{p+1} + 2\varepsilon )  \} ] ,
\end{align*}
where $\lambda_{2m-1} = \lambda_{1} $ and
$\lambda_{2m} = \lambda_{2} ,$ respectively for all integers
$m\ge 1.$ 
Then by the mean value theorem we obtain
$$ |Q| \ge |I|\cdot\exp [ -k_1 \{ 
q_{p} (\lambda_{p} + 2\varepsilon ) +
s (\lambda_{p+1} + 2\varepsilon ) \} ] .$$
Let $$\varGamma := 
\bigcap_{l=1}^{\infty}\bigsqcup_{Q\in\mathcal{Q} (l)} Q .$$
Then we have 
\begin{align*}
|S_{n_{2m-1}} \varphi (x) & - n_{2m-1}\int\varphi d\mu_1  | \\
\le &|S_{n_{2m-2}} \varphi (x)  - n_{2m-2}\int\varphi d\mu_1  | \\
&+ |S_{(q_{2m-1}+1)k_1} \varphi (f^{n_{2m-2}}(x)) 
- (q_{2m-1}+1)k_1\int\varphi d\mu_1  | \\
\le &2(n_{2m-2} +k_1)\sup_{y\in I} |\varphi (y)| 
+  q_{2m-1}k_1\beta  /8 \\
\le  &q_{2m-1} k_1 \beta  /4 \quad \le  n_{2m-1}\beta /4 ,
\end{align*}
and similarly,
$$| S_{n_{2m}} \varphi (x)  -
n_{2m}\int\varphi d\mu_2  | \le  n_{2m}\beta /4$$
for all  $x\in \varGamma$ and integers $m\ge 1 .$
Therefore,
\begin{align*}
\limsup_{n\to\infty}\frac{1}{n}S_n\varphi (x)
&-\liminf_{n\to\infty}\frac{1}{n}S_n\varphi (x) \\
\ge &| \int\varphi d\mu_1 - \int\varphi d\mu_2 | \\
&-\limsup_{m\to\infty}
| \frac{1}{n_{2m-1}} S_{n_{2m-1}} \varphi (x)  -
\int\varphi d\mu_1  | \\
&-\limsup_{m\to\infty}
| \frac{1}{n_{2m}} S_{n_{2m}} \varphi (x)  -
\int\varphi d\mu_2  | \\
\ge &\beta -\beta /4 -\beta /4 
\quad =\beta /2
\end{align*}
for all $x\in\varGamma .$
This implies that  $\varGamma\subset I_{\varphi} .$
To give a lower estimate of the Hausdorff dimension of
$\varGamma$ 
take a sequence of the probability measures 
$\nu_l :=\sum_{Q\in\mathcal{Q} (l)} \delta_{x_Q} /\sharp\mathcal{Q}(l) $ 
$(l=1,2,\ldots )$ such that
$x_Q\in Q\cap\varGamma$ holds for each $Q\in\mathcal{Q} (l)$,
and an accumulation point $\nu_0\in\mathcal{M}$ of the sequence.
Then in a similar way as in the proof of Proposition 7 
it can be checked that
\begin{align*}
\liminf_{r\to 0} &\frac{\log \nu_0 ([x-r, x+r])}{\log r}  \\
&\ge \lim_{p\to\infty}        
\min_{1\le s\le q_{p+1}}
\frac{ k_1 \{ q_{p} (h_{p}-\varepsilon )
+s(h_{p+1} -\varepsilon )\} -\log 3}
{ k_1 \{ q_{p} (\lambda_{p} + 2\varepsilon )
+ (s+1) (\lambda_{p+1} + 2\varepsilon )\} -\log |I|}   \\
&=\min\Big\{ \frac{h_1 -\varepsilon}{\lambda_1 + 2 \varepsilon }  ,  
\frac{h_2 -\varepsilon}{\lambda_2 + 2\varepsilon } \Big\}
\end{align*}
for all $x\in\varGamma$.
This implies that
\begin{equation*}
\dim_H (I_{\varphi})
\ge  \dim_H (\varGamma) 
\ge \min \Big\{ \frac{h_1 -\varepsilon}{\lambda_1 +2\varepsilon } ,
\frac{h_2 -\varepsilon}{\lambda_2 +2\varepsilon } \Big\}  
\end{equation*}
for any $\varepsilon >0.$
Letting $\varepsilon\to 0$
we obtain 
\begin{align*}
\dim_H (I_{\varphi})
\ge  \min \{ D(\mu_1 ), D(\mu_2) \} .
\end{align*}
This completes the proof of Proposition 8.

\vspace{0.2cm}

\noindent {\bf Remark.}
The existence of a return time function
satisfying (H1)-(H4) is unnecessary 
for Propositions 1,7 and 8.
The assertions of these propositions are valid 
for any topologically mixing $C^2$ map on a compact interval.


\begin{thebibliography}{99}









\bibitem{BS}{L. Barreira and J. Schmeling},
{\it Sets of "non-typical" points have 
full Hausdorff dimension and full topological entropy},
Israel J. Math. {\bf 116}
(2000), 29--70.






\bibitem{B1}{R. Bowen},
{\it Hausdorff dimension of quasi-circles},
Publ. Math. IHES {\bf 50}
(1979), 259--273.





\bibitem{BK}{H. Bruin and G. Keller},
{\it Equilibrium states for S-unimodal maps},
Ergod. Th. Dynam. Sys. {\bf 18}
(1998), 765--789.

\bibitem{C0}{ Y.M. Chung},
{\it Shadowing property of non-invertible maps with hyperbolic measures},
Tokyo J. Math. {\bf 22}
(1999),  145--166.





\bibitem{C4}{ Y.M. Chung},
{\it Recurrence times and large deviations},
preprint.




\bibitem{H}{F. Hofbauer},
{\it Local dimension for piecewise monotone maps on the interval},
Ergod. Th. Dynam. Sys. {\bf 15}
(1995), 1119--1142.

\bibitem{I}{G. Iommi},
{\it Multifractal analysis for countable Markov shifts},
Ergod. Th. Dynam. Sys. {\bf 25}
(2005), 1881--1907.




\bibitem{K}{ A. Katok},
{\it Lyapunov exponents, entropy and periodic orbits for diffeomorphisms},
Inst. Hautes Etudes Sci. Publ. Math. {\bf 51}
(1980), 137--173.



\bibitem{KM}{ A. Katok and L. Mendoza},
{\it Dynamical systems with nonuniformly hyperbolic behavior},
supplement to "Introduction to the modern theory of dynamical 
systems" written by A. Katok and B. Hasselblatt,
Cambridge Univ. Press, Cambridge, 1995, pp 659--700.







\bibitem{L2}{F. Ledrappier},
{\it Some relations between dimension and Lyapunov exponents},
Commun. Math. Phys. {\bf 81}
(1981), 229--238.





\bibitem{Na}{K. Nakaishi},
{\it Multifractal formalism for some parabolic maps},
Ergod. Th. Dynam. Sys. {\bf 20}
(2000), 843--857.




\bibitem{Ol}{L. Olsen},
{\it Multifractal analysis of divergence points of 
deformed measure theoretical Birkhoff averages},
J. Math. Pures Appl. {\bf 82}
(2003), 1591--1649.








\bibitem{P2}{Ya. Pesin},
{\it Dimension Theory in Dynamical Systems},
Univ. of Chicago Press, Chicago, 1997.

\bibitem{PW97}{Ya. Pesin and  H. Weiss},
{\it A multifractal analysis of equlibrium measures 
for conformal expanding maps and 
Moran-like geometric constructions},
J. Stat. Phys. {\bf 86}
(1997), 233--275.


\bibitem{PW}{Ya. Pesin and  H. Weiss},
{\it The multifractal analysis of Birkhoff averages 
and large deviations},
in "Global analysis of dynamical systems",
edited by H. W. Broer, B. Krauskopf and G. Vegter  
Inst. Phys., Bristol, 2001, pp 419--431.



\bibitem{PoW}{M. Pollicott and H. Weiss},
{\it Multifractal analysis of Lyapunov exponent for continued 
fraction and Manneville-Pomeau transformations and applications to
Diophantine approximation},
Commun. Math. Phys. {\bf 207}
(1999), 145--171.




\bibitem{Ru}{D. Ruelle},
{\it An inequality for the entropy  of differentiable maps},
Bol. Soc. Brasil. Math. {\bf 9}
(1978), 83--87.






\bibitem{TV2}{F. Takens and E. Verbitskiy},
{\it On the variational principle for the topological entropy
of certain non-compact sets},
Ergod Th. Dynam. Sys. {\bf 23}
(2003), 317--348.

\bibitem{W}{P. Walters},
{\it An Introduction to Ergodic Theory},
 Graduate Texts in Mathematics {\bf 79}, Springer-Verlag, New York,
1982.

\bibitem{Y0}{L.-S. Young},
{\it Dimension, entropy and Lyapunov exponents},
Ergod. Th. Dynam. Sys. {\bf 2}
(1982), 109--124.



\bibitem{Y3}{L.-S. Young},
{\it Statistical properties of dynamical systems 
with some hyperbolicity},
Ann. of Math. {\bf 147}
(1998), 585--650.







\end{thebibliography}
\end{document}